\documentclass[11pt,a4paper]{amsart}
\usepackage{graphicx,latexsym,amssymb}
\usepackage[all]{xy}
\CompileMatrices

\theoremstyle{plain}
\newtheorem{theorem*}{Theorem}
\newtheorem{theorem}{Theorem}[section]
\newtheorem{lemma}[theorem]{Lemma}

\newtheorem{proposition}[theorem]{Proposition}
\newtheorem{corollary}[theorem]{Corollary}
\newtheorem{corollary-lem}{Corollary}
\newtheorem{corollary*}[theorem*]{Corollary}
\theoremstyle{definition}

\newtheorem{example}[theorem]{Example}
\theoremstyle{remark}
\newtheorem{remark}[theorem]{Remark}

\theoremstyle{question}
\numberwithin{equation}{section}
\numberwithin{figure}{section}
\font\nb=msbm10

\def \N{\hbox{\nb N}}

\def \Z{\hbox{\nb  Z}}

\begin{document}
\title{Palindromes and orderings in Artin groups}
\author[F. Deloup]{Florian Deloup}

\address{Florian Deloup,
Laboratoire Emile Picard,
UMR 5580 CNRS/Universit\'e Paul Sabatier, 118 route de Narbonne,
31062 Toulouse, France. } \email{deloup@picard.ups-tlse.fr}

\begin{abstract}
The braid group $B_{n}$, endowed with Artin's presentation, admits
two distinguished involutions. One is the anti-automorphism
${\rm{rev}}: B_{n} \to B_{n}$, $v \mapsto \bar{v}$, defined by
reading braids in the reverse order (from right to left instead of
left to right). Another one is the conjugation $\tau:x \mapsto
\Delta^{-1}x \Delta$ by the generalized half-twist (Garside
element).

More generally, the involution ${\rm{rev}}$ is defined for all
Artin groups (equipped with Artin's presentation) and the
involution $\tau$ is defined for all Artin groups of finite type.
A palindrome is an element invariant under rev. We classify
palindromes and palindromes invariant under $\tau$ in Artin groups
of finite type. The tools are elementary rewriting and the
construction of explicit left-orderings compatible with rev.

Finally, we discuss generalizations to Artin groups of infinite
type and Garside groups.

\end{abstract}
\maketitle \tableofcontents

\section{Introduction}

\subsection{Palindromes in Artin groups}
Let $n \geq 2$. The free group $F_{n}$ on $n$ generators $s_1,
\ldots, s_{n}$ supports the involution rev$: w \mapsto
\overline{w}$ defined by
$$ s_{i_{1}}^{\alpha_{1}} \cdots s_{i_{r}}^{\alpha_{r}}
\mapsto s_{i_{r}}^{\alpha_{r}} \cdots s_{i_{1}}^{\alpha_{1}},$$
which consists in reversing the reading of the word $w$ with
respect to the prescribed set of generators. Any group $G$
presented by generators and relations which are rev-invariant
 admits such an anti-automorphism. The induced
involution will still be denoted by rev.
 The elements of $G$ which are
rev-invariant are called {\emph{palindromes}}. This
paper studies palindromes in the class of Artin groups.

A {\emph{Coxeter matrix}} of rank $n$ is a square symmetric matrix
$M$ of size $n$ with integer entries $m_{ij} \in \N \cup \{ \infty
\}$ such that $m_{ii} = 1$ for all $1 \leq i \leq n$ and $m_{ij} =
m_{ji} \geq 2$ for all $1 \leq i \not= j \leq n$. Given two
generators $a$ and $b$ of $F_{n}$ and $k \geq 2$, denote by
$w_{k}(a,b)$ the word of length $k$ defined recursively by
$$w_{2}(a,b) = ab, \ w_{k}(a,b) = \left\{
\begin{array}{cl}
w_{k-1}(a,b)b & \hbox{if}\ k\ {\hbox{is}}\ {\hbox{even}}\\
w_{k-1}(a,b)a & \hbox{if}\ k\ {\hbox{is}}\ {\hbox{odd}}.
\end{array} \right.$$
 Given a Coxeter matrix $M$, the Artin group $A_{M}$ of type $M$ is
 the group defined by the
presentation
\begin{equation} A(M) = \langle s_1, \ldots, s_n \ | \
w_{m_{ij}}(s_i, s_j) = w_{m_{ji}}(s_j, s_i) \ {\rm{for}}\
{\rm{all}}\ i\not=j \ {\rm{and}}\ m_{ij}\not= \infty. \rangle.
\label{eq:artin-pres}
\end{equation}
A group equipped with the presentation $(\ref{eq:artin-pres})$
will be called an {\emph{Artin system}} of type $M$. The set $S =
\{ s_1, \ldots, s_n \}$ is the set of positive {\emph{Artin generators}}.
 Clearly, rev fixes the word $w_{k}(a,b)$ if $k$ is even and sends it to
$w_{k}(b,a)$ if $k$ is odd. It follows that all Artin systems
carry an involutive anti-automorphism rev$:x \mapsto {\rm{rev}}(v)
= \bar{x}$ induced from the involution rev on $F_{n}$.
Accordingly, elements of the Artin group invariant under rev will
be called {\emph{palindromes}}.

Given a Coxeter matrix $M$, one can similarly define the
{\emph{Coxeter group}} $W_{M}$ of type $M$ as the quotient of the
Artin group of type $M$ by the subgroup normally generated by the
relations $s_{i}^{2} = e$, $1 \leq i \leq n$. The kernel of the
natural epimorphism $A_{M} \to W_{M}$ sending the generator $s_{i}
\in A_{M}$ to the corresponding generator $s_{i} \in W_{M}$ is
called the {\emph{pure Artin group}}.

\begin{figure}[htb]
\begin{center}
\includegraphics[height=12cm,width=8cm,angle=0,draft=false]
{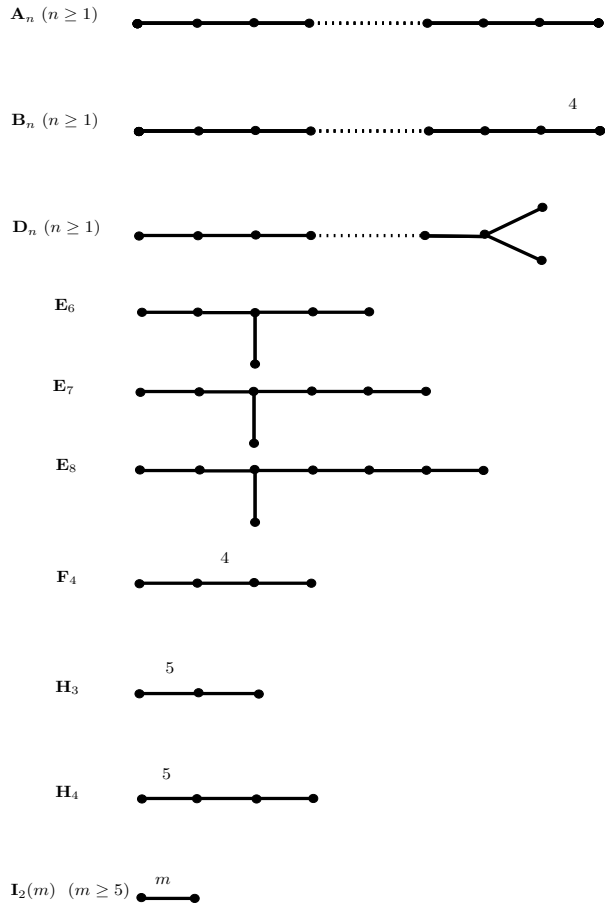} \caption{Coxeter diagrams of finite type.}
\label{fig:Coxeter}
\end{center}
\end{figure}

It is traditional to encode Artin systems in the form of a diagram
(see Figure \ref{fig:Coxeter}). Given a Coxeter matrix $M$, a
{\emph{Coxeter diagram}} is a graph $\Gamma$ whose set of vertices
is $\{ s_{1}, \ldots, s_{n} \}$ such that two vertices $s_{i},
s_{j}$ are joined by an edge if and only if $m_{ij} \geq 3$ and an
edge between two vertices $s_{i}, s_{j}$ is labelled by $m_{ij}$.
(It is customary to omit the label if $m_{ij} = 3$.) We shall
index Artin and Coxeter groups indifferently with the Coxeter
diagram $\Gamma$ or with the Coxeter matrix $M$.
 Since the Artin and Coxeter groups associated to disjoint
Coxeter diagrams are direct products, we shall assume in this paper
that all Coxeter diagrams are connected.  An Artin group
$A_{\Gamma}$ (resp. a Coxeter diagram $\Gamma$) is said to be of
{\emph{finite type}} if the associated Coxeter group $W_{\Gamma}$
is finite.

As an example, $A_{{\bf{A}}_{n}} = B_{n+1}$, the classical braid
group on $n+1$ strands. In this case, the associated Coxeter group
is the symmetric group $S_{n+1}$.

Let $A_{\Gamma}$ be an Artin system with set $S$ of Artin generators.
The Artin monoid $A_{\Gamma}^{+}$ is the
{\emph{monoid}} presented by the same generators and relations as
the Artin system $A_{\Gamma}$. Let $I$ be a subset of $S$.
Suppose that all generators in $I$ have a (right)
common multiple. Then there is a uniquely determined least common
multiple, called the {\emph{fundamental element}} $\Delta_{I} \in
A_{\Gamma}^{+}$ for the subset $I$. Such a fundamental element is
always palindromic. A basic result in the theory of Artin monoids
asserts that the fundamental element $\Delta_{S}$ exists for the
whole set $S = \{ s_1, \ldots, s_n \}$ of Artin generators if and
only if $\Gamma$ is of finite-type. It is known that $\Delta_{S}$
or $\Delta^{2}_{S}$ lies in the center of $A_{\Gamma}$. In
particular, the automorphism
\begin{equation}
\tau: A_{\Gamma} \to A_{\Gamma}, \  \ x \mapsto \Delta_{S}^{-1}x
\Delta_{S} \label{eq:def-tau}
\end{equation}
is always of order at most two. In this paper, we study and give
classification results for palindromes in Artin groups and for
$\tau$-invariant palindromes in Artin groups of finite type. We
also indicate generalizations to Artin groups of infinite type and
Garside groups.

\begin{figure}[htb]
\begin{center}
\includegraphics[height=6cm,width=8cm,angle=0,draft=false]{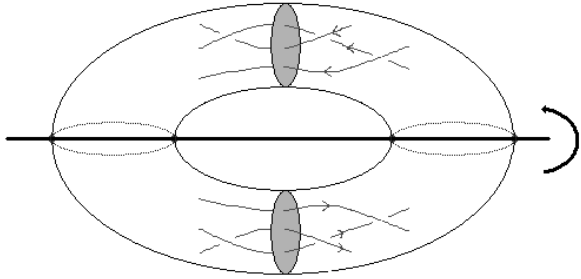}
\caption{The Weierstrass involution and palindromic braids.}
\label{fig:hyperelliptic}
\end{center}
\end{figure}

Palindromes and $\tau$-invariant palindromes have nice geometric
interpretations for $A_{\Gamma} = B_{n}$ (the classical braid
group on $n$ strands). Given a geometric braid $x$, denote by
$\widehat{x}$ its closure into a link inside a fixed solid torus
$D^{2} \times S^{1}$. The solid torus admits the involution
$$ {\rm{inv}}: D^{2} \times S^{1} \to D^{2} \times S^{1}, \ \ (re^{it},
\theta) \mapsto (r e^{-it}, - \theta),$$ whose set of fixed points
consists of two segments ($t \equiv 0$ (mod $\pi$) and $\theta
\equiv  0$ (mod $\pi$)), which is the intersection of the axis of
the 180$^{\rm{o}}$ rotation with the solid torus. Restricted to
the boundary, this is just the Weierstrass involution of the
standard torus. Observe that $\widehat{{\rm{rev}}(x)}$ is nothing
else than inv$(\widehat{x})$ with the opposite orientation. In
particular, if a braid $x \in B_{n}$ is palindromic then
$\widehat{x}$ coincides with ${\rm{inv}}(\widehat{x})$ with the
opposite orientation, see Figure \ref{fig:hyperelliptic}.

\begin{figure}[!htb]
\begin{center}
\includegraphics[height=7cm,width=9cm,angle=0,draft=false]{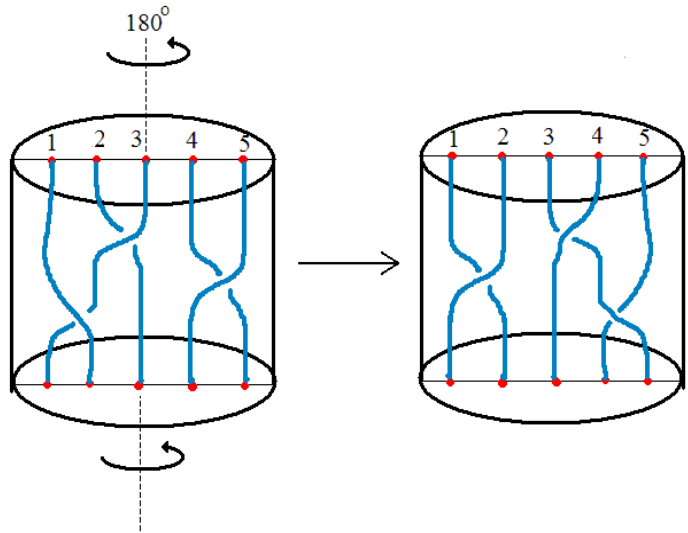}
\caption{The involution $\tau$ is induced by a vertical
$180^{\rm{o}}$ rotation of the cylinder.} \label{fig:tau}
\end{center}
\end{figure}

The fundamental element coincides with the Garside element $\Delta
= \Delta_n \in B_{n}$, the generalized half-twist on $n$ strands,
defined inductively by
$$ \Delta_2 = s_1, \ {\rm{and}}\
\Delta_n = s_1 s_2 \cdots s_{n-1}\cdot \Delta_{n-1},$$ where
$s_{1}, \ldots, s_{n-1}$ are the Artin generators of $B_n$. It
turns out, as is directly verified, that $\tau(s_{i}) = s_{n-i}$
for $1 \leq i \leq n-1$. Thus any braid commutes with $\Delta^{2}$
(a full twist); $\tau$-invariant braids are those which commute
with $\Delta$. See Figure \ref{fig:tau}.

\subsection{Statement of results}
Let $A_{\Gamma}$ be an Artin group. There is a canonical way of
producing palindromes from $A_{\Gamma}$. It consists in applying
the map
$$ {\mathfrak{pal}}: A_{\Gamma} \to A_{\Gamma}, \ \ x \mapsto x \cdot {\overline{x}} $$
which we call the {\emph{palindromization}} map.

\begin{theorem} \label{th:palindromization}
Let $\Gamma$ be of finite type. The palindromization map
${\mathfrak{pal}}:A_{\Gamma} \to A_{\Gamma}$ is injective and its
image is the subset of pure palindromes.
\end{theorem}

The injectivity of ${\mathfrak{pal}}$ was first proved in
\cite{DGKT} in the case of the classical braid group $B_{n}$,
using the Jacquemard algorithm. \\

The proof of Theorem \ref{th:palindromization} is given in \S
\ref{sec:pure}. Here is an outline. The proof is partially based
on the existence of a left-invariant ordering of a certain type on
Artin groups (Theorem \ref{th:left-orders} below). Recall that a
group $G$ is {\emph{left-ordered}} (or has a left-invariant
ordering) if there exists a total ordering $<$ on the set $G$ such
that $x < y$ implies $ax < ay$ for all $x, y, a \in G$.

\begin{lemma} \label{lem:general}
Let $G$ be a left-ordered group, equipped with an
involutive anti-automorphism $G
\to G, x \mapsto \overline{x}$, such that $x > e_{G}$ if and only
if $\overline{x} > e_{G}$. Then the palindromization map $x
\mapsto x \overline{x}$ is injective.
\end{lemma}

\begin{remark}
Actually, Lemma \ref{lem:general} is true for any automorphism or
anti-automorphism $x \mapsto \overline{x}$.
\end{remark}

We show the existence of an explicit left-invariant ordering on
Artin groups of type ${\bf{A}}$, ${\bf{B}}$ and ${\bf{D}}$,
``compatible" with the involution rev:

\begin{theorem} \label{th:left-orders}
Let $\Gamma = {\bf{A}}_n, {\bf{B}}_n$ or ${\bf{D}}_n$. There
exists a left-invariant ordering $<$ on $A_{\Gamma}$ such that $x
> e$ if and only if $\overline{x} > e$ for all $x \in A_{\Gamma}$.
\end{theorem}

In fact, there exists such an ordering on $B_n$, namely the
Dehornoy ordering. (Thurston-type orderings all have that
property.) We use it to prove Theorem \ref{th:left-orders}. Then
we embed other Artin groups of finite type into these in such a
way that we can apply Lemma \ref{lem:general}. We do not know
about the left-orderability of the Artin group associated to
${\bf{E}}_{8}$. Hence we cannot conclude at the present time for
some Artin groups of finite type by this method. Hence we resort
to the Jacquemard algorithm and the combinatorial method
developped in \cite{DGKT} to prove that Theorem
\ref{th:palindromization} is true for $\Gamma = {\bf{E}}_{8}$.
Then we embed all remaining Artin groups of finite type into
$A_{{\bf{E}}_{8}}$ and finish the proof.

\begin{remark}
We also give in \S \ref{sec:remarks} a proof outline valid
for all Artin groups of infinite type using the Jacquemard
algorithm.
\end{remark}

\begin{remark} \label{rem:notall}
It is readily seen that not all palindromes are in the image of
${\mathfrak{pal}}$, even for the classical braid group. For
instance, $\Delta \not\in {\rm{Im}}({\mathfrak{pal}})$ (in fact,
$\Delta$ is not even a pure braid) and ${\overline{\Delta}} =
\Delta$. Fig. \ref{fig:Delta}
 below displays the equality
$\Delta = {\overline{\Delta}}$ for $\Delta = \sigma_{1} \sigma_{2}
\sigma_{3} \sigma_{1} \sigma_{2} \sigma_{1}$, the generalized
half-twist on four strands. More generally, for all Artin groups
of finite type, $\Delta_{S}$ is square-free (it is not represented
by a word which contains a square $s^{2}$, $s \in S$), so
$\Delta_{S}$ is never in the image of ${\mathfrak{pal}}$.
\end{remark}

\begin{figure}[htb]
\begin{center}
\includegraphics[height=1.5cm,width=8cm,angle=0,draft=false]
{Delta.pstex} \caption{A palindromic braid not of the form
$v\overline{v}$.} \label{fig:Delta}
\end{center}
\end{figure}

We note the following consequences of Theorem
\ref{th:palindromization}.

\begin{corollary} \label{lem:delta-associated}
Let $x \in A_{\Gamma}$ such that ${\rm{rev}} \circ \tau (x) = x$.
Assume that $x$ has the same image in $W_{\Gamma}$ as $\Delta =
\Delta_{S}$. Then there exists $\delta \in A_{\Gamma}$ such that
$x = \Delta \delta\overline{\delta} =
\tau(\delta)\Delta\overline{\delta}$.
\end{corollary}

\noindent{\bf{Proof}}. Observe that $y = \Delta^{-1} x$ is pure
and
$$\overline{y}  =  \overline{x}\overline{\Delta}^{-1}  =
\overline{x} \Delta^{-1}  =  \tau(x) \Delta^{-1}
     =  \Delta^{-1} x
     =  y.$$
We can therefore apply Theorem \ref{th:palindromization}: there
exists $\delta \in A_{\Gamma}$ such that $y = \Delta^{-1} x  =
\delta
\overline{\delta}$. \hfill $\blacksquare$\\

\begin{corollary} \label{cor:pure-rev-tau}
Let $x \in A_{\Gamma}$ pure and invariant under ${\rm{rev}}$ and
$\tau$. Then there exists a unique $\delta \in A_{\Gamma}$ such
that $x = \delta \overline{\delta}$ and $\tau(\delta) = \delta$.
\end{corollary}

\noindent{\bf{Proof}}. The existence of $\delta \in A_{\Gamma}$
such that $x = \delta \overline{\delta}$ follows from Theorem
\ref{th:palindromization}. Since ${\rm{rev}}$ and $\tau$ commute,
we have $\tau(x) = \tau(\delta) \overline{\tau(\delta)} = \delta
\overline{\delta} = x$. Applying again Theorem
\ref{th:palindromization}
 (injectivity of ${\mathfrak{pal}}$) yields
$\tau(\delta) = \delta$. \hfill $\blacksquare$\\

We now describe a general decomposition for palindromes in an
Artin group $A_{\Gamma}$ of finite type. Recall that $S$ denotes
the set of positive Artin generators.


\begin{theorem} \label{th:dec-pal}
Let $x \in A_{\Gamma}$ be a palindrome. Then there exist $y \in
A_{\Gamma}$ and $I \subseteq S$ such that
\begin{equation}
x =  y \ \Delta_{I} \ \overline{y}. \label{eq:dec-pal}
\end{equation}
Suppose that $A_{\Gamma}$ is left-ordered with $<$ and the
restriction of the order to $A_{\Gamma}^{+}$ is a well-ordering.
Then the decomposition $(\ref{eq:dec-pal})$ is unique provided
that
 $(\Delta_{I}, y)$ is minimal with respect to the lexicographic ordering
$\prec\ = ( < , <)$ on $A_{\Gamma} \times A_{\Gamma}$.
\end{theorem}

The decomposition (\ref{eq:dec-pal}) with the foregoing
requirement on $(\Delta_{I}, \gamma)$ will be called the
{\emph{canonical}} decomposition.
If we denote by $<^{\rm{opp}}$
the oppositive ordering on $A_{\Gamma}$ ($ x <^{\rm{opp}} y$
if and only  if $y < x$), then there is also a smilar
decomposition with the requirement that $(\Delta_I, y)$
be minimal with respect to the ordering $( <^{\rm{opp}} , <)$. \\

\begin{remark}
The main result of \cite{Laver} asserts that there exists a
left-ordering (the Dehornoy ordering, see \S \ref{sec:pure} for a
definition) on $B_{n}$ whose restriction to $B_{n}^{+}$ is a
well-ordering. Any left-ordering of Thurston type has also this
property \cite[Prop. 7.3.1]{DSW}. This result can be extended to
Artin groups $A_{\Gamma}$ for $\Gamma = {\bf{B}}_{n}$ or
${\bf{D}}_{n}$: there is a left-ordering on $A_{\Gamma}$ whose
restriction to $A_{\Gamma}^{+}$ is a well-ordering. I do not know
if this is true for other Artin groups.
\end{remark}

\begin{remark}
Theorem \ref{th:dec-pal} yields a partial ordering on the set of
palindromes. It does not coincide with the initial left-invariant
ordering restricted to the subset of palindromes. In the case when
$A_{\Gamma} = B_{n}$, if the ordering on $B_{n}$ is the Dehornoy
ordering, then the ordering of elements $(e,y) \prec (\Delta_{I},
y) \prec (\Delta_J, y)$ does coincides with the Dehornoy ordering
restricted to the corresponding palindromes: $y \bar{y} < y
\Delta_I \bar{y} < y \Delta_J \bar{y}$ for $\Delta_I < \Delta_J$.
(Note that if we identify $I, J \subseteq S = \{ s_1, \ldots,
s_{n-1} \}$ to subsets of $\{ 1, \ldots, n-1 \}$, then $\Delta_I <
\Delta_J$ is equivalent to $I > J$ in the usual lexicographic
ordering of subsets of $\{1, \ldots, n-1\}$.)
\end{remark}

\begin{remark}
Suppose that $\alpha \in B_{n}^{+}$ is palindromic. To ensure
uniqueness of the decomposition $(\ref{eq:dec-pal})$ for $\alpha$,
it is not enough to require that the length $\ell(y)$ of $y$ be
extremal, as the following example shows:
$$ \alpha = (\sigma_3 \sigma_5)\Delta_{\{1, 2\}}(\sigma_5 \sigma_3)
= (\sigma_5 \sigma_1)\Delta_{\{2,3\}}(\sigma_1 \sigma_5) \in
B_6.$$ In fact, it is not even enough to fixe $\Delta_{I}$ in the
decomposition to ensure uniqueness, as the following example
shows.
\end{remark}

\begin{example} \label{ex:fondamental}
Consider the braid $x = (\sigma_2 \sigma_3 \sigma_1)^{2} \in
B_{4}^+$. It is readily verified that $x = \Delta_{\{1, 2, 3\}}$,
hence $x$ is palindromic. However, this is not the canonical
decomposition. The following equality shows that the map $y
\mapsto y \Delta_{I} \overline{y}$ is not injective in general:
$$ \begin{array}{cccccc}
x =  \sigma_3 \sigma_2 \sigma_1 \sigma_3 \sigma_2 \sigma_3
 & = &  \sigma_2 \sigma_3 \sigma_2 \sigma_1
\sigma_2 \sigma_3  & = &  \sigma_2 \sigma_3 \sigma_1
\sigma_2 \sigma_1 \sigma_3  \\
& & & = &  \sigma_2 \sigma_1 \sigma_3 \sigma_2 \sigma_3
\sigma_1  \\
& &  & = &  \sigma_2 \sigma_1 \sigma_2
\sigma_3 \sigma_2 \sigma_1  \\
& & & = &  \sigma_1 \sigma_2 \sigma_1 \sigma_3 \sigma_2 \sigma_1
\end{array}$$
Therefore, $x =  \sigma_3 \sigma_2 \Delta_{\{1, 3\}} \sigma_2
\sigma_3 =  \sigma_1 \sigma_2 \Delta_{\{1, 3\}} \sigma_2
\sigma_1.$ If we endow $B_{n}$ with the Dehornoy ordering, then
$\sigma_3 \sigma_2 < \sigma_{1} \sigma_{2}$, so the second
decomposition is not the canonical decomposition. We leave it to
the reader to verify that the first decomposition is the canonical
decomposition.
\end{example}

\section{Pure palindromes} \label{sec:pure}

Although it is possible to present a slightly more direct proof of
Theorem \ref{th:palindromization}, the argument we present
establishes stronger results about orderings on Artin groups of finite type.\\

\noindent{\bf{Proof of Lemma \ref{lem:general}}}. Assume for
instance that $x
> e_{G}$. By assumption, $\bar{x} > e_{G}$. Then $x \bar{x} >
e_{G}$. We have just proved that $x \bar{x} = e_{G}$ implies $x =
e_{G}$. Assume now that $x \bar{x} = y \bar{y}$. Write $y = x z$,
$z \in A_{\Gamma}$. Then we have $\overline{y} = \overline{z}
\overline{x}$, hence $y \overline{y} = x\ \gamma\ \overline{z}\
\overline{x} = x \overline{x}$. Therefore $z \overline{z} = e_G$.
By our previous argument, $z = e_G$. Hence ${\mathfrak{pal}}$ is
one-to-one. \hfill $\blacksquare$\\

We now turn to the proof of Theorem \ref{th:palindromization}
proper. We focus in this section on the proof of the first
statement (injectivity of ${\mathfrak{pal}}$). The proof of the
second statement (that the image of ${\mathfrak{pal}}$ coincides
with pure palindromes) is postponed after the proof of Theorem
\ref{th:dec-pal}.\\

\noindent{\bf{Step 1}} (Proof of Theorem \ref{th:palindromization}
for the braid group $B_n$). There is a total left-ordering $<$ on
$B_{n}$, called the Dehornoy order. A word of the form
$$x_{0}s_i x_{1}s_i \cdots x_{k-1} s_{i}
x_{k}$$ where $x_{1}, \ldots, x_{k}$ are words in the letters
$s_{i+1}^{\pm 1}, \ldots, s_{n-1}^{\pm 1}$, is called a
$s_i$-{\emph{positive}} word. A braid $x$ is
$s_i$-{\emph{positive}} if $x$ can be represented by a
$s_i$-positive word. A braid is called $s_{i}$-negative if its
inverse is $s_{i}$-positive. Call a braid $s$-{\emph{positive}} or
$s$-{\emph{negative}} if it is $s_i$-positive or $s_{i}$-negative
for some $i$.

\begin{theorem}[Dehornoy, \cite{Dehornoy}] \label{th:Dehornoy}
There is a set partition of $B_{n}$ into three classes:
$s$-positive braids, $s$-negative braids, and the trivial braid.
\end{theorem}

\noindent The (left) Dehornoy order is defined by setting $x < y$
if and only if $x^{-1}y$ is $s$-positive. It is obvious that
$B_{n}$ endowed with the Dehornoy ordering satisfies the condition
of Lemma \ref{lem:general}. Therefore, Lemma \ref{lem:general}
applies. \hfill $\blacksquare$\\

\noindent{\bf{Step 2}} (Construction of left-orderings and proof
of Theorem \ref{th:palindromization} for Artin groups of type
${\bf{B}}$ and ${\bf{D}}$). Let $(G,<~)$ be a left-invariant
ordered group and let $\varphi$ be an automorphism or an
anti-automorphism of $G$. A pair $(\varphi, <)$ will be said to
have {\emph{property}} (PPC)\footnote{preserves the positive
cone.} (resp. (SPPC)\footnote{strongly preserves the positive
cone.}) when $\varphi(x)> e$ if (resp. and only if) $x > e$ for
all $x \in G$. Our first goal is to construct a pair satisfying
(SPPC), so as to apply again Lemma \ref{lem:general}.

\begin{lemma} \label{lem:sppc}
Let $G$ be a group equipped with an automorphism or
anti-automorphism $\varphi$, sitting in a short exact sequence
$$\xymatrix{
1 \ar[r] & H \ar[r]^{i} & G \ar[r]^{p} & K \ar[r] & 1}$$ of groups
such that $(H,<)$ and $(K,<)$ are both left-invariant ordered
groups and $\varphi(i(H)) \subseteq i(H)$. Suppose that the
automorphism or anti-automorphism $\varphi_{H} = i^{-1} \circ
\varphi \circ i$ has the property {\rm{(PPC)}} (resp.
{\rm{(SPPC)}} and that $p \circ \varphi(x) > e$ if and only if
$p(x) > e$ for all $x \in G$. Then there is a left-invariant order
$<$ on $G$ such that $(\varphi, <)$ is {\rm{(PPC)}} (resp.
{\rm{(SPPC)}} pair.
\end{lemma}

\noindent{\bf{Proof}}. Let $x, y \in G$. We declare $x < y$ if
$p(x) < p(y)$ or else $p(x) = p(y)$ and $e < i^{-1}(x^{-1}y)$.
This defines a left-invariant order on $G$. The claimed properties
are left to the reader to verify. \hfill $\blacksquare$\\

\begin{remark}
If $\varphi$ is periodic, that is $\varphi^{k} = {\rm{id}}$ for
some $k \geq 1$, then (PPC) is equivalent to (SPPC). This is the
case in particular for $\varphi = {\rm{rev}}$.
\end{remark}

Let $G$ be an Artin system of type ${\bf{B}}$ or ${\bf{D}}$. There
is a natural projection ${\pi}:G \to B_{n}$, easily described in
terms of the Artin generators. For convenience, denote by
$\beta_1, \ldots, \beta_n$ (resp. $\delta_1, \ldots, \delta_n$)
the Artin generators of the Artin system $A_{{\bf{B}}_{n}}$ (resp.
$A_{{\bf{D}}_{n}}$). We keep the notation $s_1, \ldots, s_{n-1}$
for the generators of the braid group $B_{n}$.

If $G = A_{{\bf{B}}_{n}}$, then $$\pi(\beta_j) = \left\{
\begin{array}{cl} \sigma_j & {\rm{if}}\ 1 \leq j \leq n-1 \\
e & {\rm{if}}\ j = n \end{array}\right..$$

If $G = A_{{\bf{D}}_{n}}$ then $$\pi(\delta_j) = \left\{
\begin{array}{cl} \sigma_j & {\rm{if}}\ 1 \leq j \leq n-2 \\
\sigma_{n-1} & {\rm{if}}\ j = n-1,\ n \end{array}\right..$$

Therefore there is a short exact sequence
$$ \xymatrix{
1 \ar[r] & {\rm{Ker}}(\pi) \ar[r] & G \ar[r]^{\pi} & B_{n} \ar[r]&
1. }$$ It is shown in \cite{CrispParis} that Ker$(\pi)$ is a free
group of rank $n$ or $n-1$ (according to whether $G =
A_{{\bf{B}}_{n}}$ or $G = A_{{\bf{D}}_{n}}$).

\begin{lemma} \label{lem:rev-free}
Let $F_{n}$ be a free group of order $n$. There is a
left-invariant order $<$ on $F_{n}$ such that $({\rm{rev}}, <)$
has $({\rm{SPPC}})$.
\end{lemma}

\noindent{\bf{Proof}}. We use the Magnus ordering defined as
follows. Let $F_{n}$ be freely generated by $x_1, \ldots, x_{n}$
and $\Lambda = \Z \langle \langle X_1, \ldots, X_n \rangle
\rangle$ be the ring of formal power series in the non-commuting
indeterminates $X_1, \ldots, X_{n}$. The Magnus map $\mu:F_{n} \to
\Lambda$ defined by $$\mu(x_{j}) = 1 + X_{j}, \ \mu(x_{j}^{-1}) =
1 - X_{i} + X_{i}^{2} - X_{i}^{3} + \cdots$$ is an injective
mapping of $F_{n}$ into the multiplicative subgroup of series
whose first coefficient (degree $0$ coefficient) is $1$. Now order
$\Lambda$ as follows. We first list formal power series according
to the total degree of monomials. Now monomials $X_{j_{1}}^{k_{1}}
\ldots X_{j_{r}}^{k_{r}}$ of a given degree $d = \sum_{1 \leq i
\leq r} k_{i}$ are ordered lexicographically according to the
r-uple of subscripts $(j_1, \ldots, j_r)$. Then two series are
compared by looking at the first term at which their coefficient
differ and order them according to that coefficient. This defines
a total order on $\Lambda$ whose restriction to the image of $\mu$
is left-invariant (in fact even bi-invariant). There is a kwown
sufficient condition for the Magnus ordering to have (SPPC).

\begin{lemma}\label{lem:magnus-sppc}
Let $\varphi:F_{n} \to F_{n}$ be an automorphism or an
auti-automorphism. If the induced map $\varphi_{{\rm{ab}}}:
F/[F_{n}, F_{n}] \to F/[F_n, F_n]$ on the abelianization of $F_n$
is the identity, then the Magnus ordering has {\rm{(SPPC)}}.
\end{lemma}

A proof can be found in \cite{DSW}, Proposition 9.2.5. The proof
there is given for an automorphism $\varphi$ but works as well if
$\varphi$ is an anti-automorphism. Since rev clearly satisfies the
hypotheses of Lemma \ref{lem:magnus-sppc}, application of Lemma
\ref{lem:magnus-sppc} concludes the proof of Lemma
\ref{lem:rev-free}. \hfill $\blacksquare$\\

\noindent We shall apply Lemma \ref{lem:sppc} to $(G,
{\rm{rev}})$. Clearly, the hypotheses pertaining to $H = F_{n}$
(or $F_{n-1}$), equipped with the Magnus ordering, and to $K =
B_{n}$, equipped with the Dehornoy ordering, are satisfied.
Applying Lemma \ref{lem:sppc}, we obtain a left-ordering on $G$
with respect to which rev has (SPPC). Therefore, we can apply
Lemma \ref{lem:general} to derive
the desired conclusion. \hfill $\blacksquare$\\

\noindent{\bf{\underline{Step 3}}} (Embeddings of Artin groups).
One general method to construct a left-ordering on a group $G$ is
to embed it into a left-orderable group. For our purpose, it is
sufficient to construct an embedding with a special property.

\begin{lemma} \label{lem:special-embedding}
Let $i:(H, \varphi_{H}) \hookrightarrow (G, \varphi)$ be an
embedding of groups equipped with anti-automorphisms such that $i
\circ \varphi_{H} = \varphi \circ i$. If $G$ is left-ordered and
has {\rm{(SPPC)}}, then $H$ has the same properties.
\end{lemma}

It is known that the Artin groups of type {\bf{H}}$_{\mathbf{3}}$
and {\bf{I}}$_{\mathbf{2}}$ inject into Artin groups of type
${\bf{D}}$ and ${\bf{A}}$ respectively. Furthermore, the
embeddings can be realized \cite{Crisp} so as to ensure that the
images of the Artin generators are palindromes (invariant under
rev). Therefore, since Artin groups of type ${\bf{D}}$ and
${\bf{A}}$ equipped with rev have (SPPC), Lemma
\ref{lem:special-embedding} applies for Artin groups associated to
{\bf{H}}$_{\mathbf{3}}$ and {\bf{I}}$_{\mathbf{2}}$ with $\varphi
= {\rm{rev}}$. Then Lemma {\ref{lem:general}} applies. \hfill
$\blacksquare$\\

\noindent{\bf{\underline{Step 4}}} (Other Artin groups of finite
type). If we knew that the Artin group associated to
${\bf{E}}_{8}$ is left-orderable with property (SPPC), then we
could apply the previous argument (step 3) to all remaining Artin
groups of finite type. Unfortunately, we do not know whether this
is true. So we use a different method for the remaining cases. It
is based on the following observation.

\begin{lemma} \label{lem:embeddings}
Let $i:A_{\Gamma'} \hookrightarrow A_{\Gamma}$ be an embedding of
Artin groups such that $i \circ {\rm{rev}} = {\rm{rev}} \circ i$.
If the palindromization map ${\mathfrak{pal}}: A_{\Gamma} \to
A_{\Gamma}$ is one-to-one, then the palindromization map
${\mathfrak{pal}}: A_{\Gamma'} \to A_{\Gamma'}$ is also
one-to-one.
\end{lemma}

The remaining Artin groups (for which the previous steps do not
apply) are associated to Coxeter diagrams ${\bf{E}}_{6},
{\bf{E}}_{7}, {\bf{E}}_{8}, {\bf{F}}_{4}$ and ${\bf{H}}_{4}$. All
these groups admit an embedding into $A_{{\bf{E}}_{8}}$ such that
the images of Artin generators are palindromes. By Lemma
\ref{lem:embeddings}, it is therefore sufficient to prove that
${\mathfrak{pal}}:A_{{\bf{E}}_{8}} \to A_{{\bf{E}}_{8}}$ is
one-to-one. Then by the usual argument (briefly recalled below),
since the natural map $A_{\bf{{E}}_8}^{+} \to A_{{\bf{E}}_{8}}$ is
an embedding, it is sufficient to prove that the restriction of
${\mathfrak{pal}}$ to $A_{{\bf{E}}_{8}}^{+}$ is one-to-one.

We recall an algorithm due to A. Jacquemard \cite{Jacquemard}. It
is originally expressed for the classical braid monoid
$B_{n}^{+}$, but all his arguments go mutatis mutandis for
$A_{{\bf{E}}_{8}}^{+}$ (and other Artin monoids). The input is a
pair $(w,s)$ where $w$ is a word representing an element in
$A_{{\bf{E}}_{8}}^{+}$ written in letters representing positive
Artin generators and $s \in S$ is one letter representing a
positive Artin generator. The algorithm decides whether there is a
new word $w' = sw"$ that represents the same element in
$A_{{\bf{E}}_{8}}^{+}$ and that starts with the letter $s$. The
corresponding Artin generator $s$ is said to be
{\emph{left-extractible}} from $w \in A_{{\bf{E}}_{8}}^{+}$. If
this is the case, it returns the new word $w'$. If not, it returns
"false". The algorithm is recursive and based on two steps: 1)
Swap $s$ with the immediate left neighbor $s'$ as long as $s$ and
$s'$ commute when regarded in $A_{{\bf{E}}_{8}}^{+}$. If $s$
becomes the first letter of the word, we are done. This step stops
if and only if $w = w_0 s' s w_1$ where $w_0$ is a non-empty word
and $s' s s' = s s' s$ in $A_{{\bf{E}}_{8}}^{+}$. 2) This step
occurs only after a call to the first step has been performed. The
algorithm calls itself recursively with the pair $(w_1, s')$. If
the call is successful, the new word is of the form $w = w_0 s' s
s' w_2$, so that we can apply the Coxeter relation to turn it into
$w = w_0 s s' s w_2$ and return to the first step. If the call
fails to rewrite $w_1 = s' w_{1}'$, then the algorithm stops and
returns false.

Define the {\emph{blocking left-index}} $I_{s}(w)$ of $s$ in the
word $s$ to be the number of letters $s'$ in the word $w$ to the
left of $s$ such that $s' s s' = s s' s$ in $A_{{\bf{E}}_8}^{+}$.
(This notion does depend on the word $w$, not on the element it
represents.) Since either the blocking left-index of $s$ strictly
decreases [step 1] or the length of the word $w$ (to which the
algorithm applies) strictly decreases [step 2], the algorithm
terminates.

The crucial point that requires to be verified is the following
claim. Let $w = w_0 s' s w_1$ be a word such that $w_0$ is a
non-empty word that does not contain the letter $s$ and where
(when viewed in $A_{{\bf{E}}_{8}}^{+}$) $s' s s' = s s' s$. The
claim is: if $s'$ is not left-extractible in $w_1$, then $s$ is
not left-extractible in $w$. Following Jacquemard's original
argument, we see that though the blocking left-index of $s$ may
decrease, it cannot be zero because after each word
transformation, there will always be a letter $s'$ occurring on
the left before the leftmost letter $s$.

We now review the key lemma (Lemma 3.2) of \cite{DGKT}. The Lemma
is proved there in the setting of the braid monoid, but the proof
continues {\emph{mutatis mutandis}} to be valid in the case of
$A_{{\bf{E}}_8}^{+}$. Suppose
that $x \bar{x} = s y \bar{y} s$ in $A_{{\bf{E}}_8}^{+}$ for some
positive Artin generator $s$. Let $w_{0} = x \bar{x}, w_{1},
\ldots, w_{k} = s y'$ be a finite sequence of positive words
representing $x \bar{x}$ such that each word $w_{i}$ is obtained
from $w_{i-1}$ by a relation in $A_{{\bf{E}}_8}^{+}$ according to
Jacquemard's algorithm. Then each relation is performed only
within the first half of the word $w_i$, $i = 0, 1, \ldots, k,$
which implies that all relations involve only letters from $x$.\\

Using this lemma, an induction on the length of the positive word
$w$ shows that $v \bar{v} = w \bar{w}$ (for a positive word $v$)
in $A_{{\bf{E}}_8}^{+}$ if and only if $v = w$ in
$A_{{\bf{E}}_8}$. (Here we use the fact that $A_{{\bf{E}}_8}^{+}$
embeds in $A_{{\bf{E}}_8}$.) The general case follows since given
any $x \in A_{{\bf{E}}_8}$, there exists a central element
$\Delta^{N} \in A_{{\bf{E}}_8}^{+}$ such that $\Delta^{N} x \in
A_{{\bf{E}}_8}^{+}$. Therefore, Theorem \ref{th:palindromization}
is true for $\Gamma = {\bf{E}}_{8}$ and finally, for the remaining
$\Gamma = {\bf{E}}_{6}, {\bf{E}}_{7}, {\bf{F}}_{4}$ and
${\bf{H}}_{4}$. This achieves the proof.  \hfill $\blacksquare$\\





\section{General palindromes}
\label{sec:classification-pal}

This section is devoted to the proof of Theorem \ref{th:dec-pal}
and the fact that the image of ${\mathfrak{pal}}$ consists of the
subset of pure palindromes.

\subsection{Proof of Theorem \ref{th:dec-pal}}

\subsubsection{Existence} The existence is the consequence of the
following two lemmas. Given $x \in A_{\Gamma}^{+}$, denote by
$S(x) = \{ s \in S \ | \ x = sy\ {\hbox{for some}}\ y \in
A_{\Gamma}^{+} \}$ (starting set) and $F(x) = S(\overline{x})$
(finishing set). The first lemma follows from the divisibility
theory for Artin monoids.

\begin{lemma} \label{lem:prefix}
Let $x \in A_{\Gamma}^{+}$. If $J \subseteq S(x)$, then there
exists $y \in A_{\Gamma}^{+}$ such that $x = \Delta_{J} y$.
\end{lemma}

\begin{lemma} \label{lem:pseudo-dicho}
Let $x \in A_{\Gamma}^{+}$ be a palindrome. Then there exists $J
\subseteq S$ such that $x = \Delta_{J}$ or there exists $s \in S$
such that $x = s a s$ for some $a \in A_{\Gamma}^{+}$.
\end{lemma}

\noindent{\bf{Proof of lemma \ref{lem:pseudo-dicho}}}. Apply Lemma
\ref{lem:prefix} to $x$ and $J = S(x)$, thus $x = \Delta_{S(x)} y$
for some $y \in A_{\Gamma}^{+}$. If $y = e$, we are done.
Otherwise $F(y) \not= \varnothing$, so there exists $s \in F(y)$.
Since $\overline{x} = x$, we have $F(y) \subseteq S(x)$. Hence $s
\in S(x)$. Thus there are $\delta, y' \in A_{\Gamma}^{+}$ such
that $s\delta = \Delta_{S(x)}$ and $y's = y$. Then $x = s
a s$ for some $a = \delta y'$. \hfill $\blacksquare$\\

Let now $x \in A_{\Gamma}^{+}$ such that $\overline{x} = x$. We
apply Lemma \ref{lem:pseudo-dicho} to $x$. If $x = \Delta_{J}$ for
some $J \subseteq S$, we are done. Otherwise $x = sa_{1}s$ for
some $a_{1} \in A_{\Gamma}^{+}$. Since $sa_{1}s = x = \overline{x}
= s \overline{a_{1}} s$, we deduce that $\overline{a_{1}}
  = a_{1}$. So Lemma \ref{lem:pseudo-dicho} applies to $a_{1}$.
  Denote by $\ell(x)$ the length (in positive Artin generators) of
  an element $x \in A_{\Gamma}^{+}$.
An immediate induction using repeatedly Lemma
\ref{lem:pseudo-dicho} shows that either $x = y \Delta_{I}
\overline{y}$ for some $y \in A_{\Gamma}^{+}$ with $\ell(y) <
\ell(x)/2$ or $x = y a \overline{y}$  where $y, a \in
A_{\Gamma}^{+}$ and $\ell(y) = \ell(x)/2$. In the latter case,
since
$$\ell(x) = \ell(y a \overline{y})=
\ell(\gamma) + \ell(a) + \ell(\overline{y}) = 2\ell(y) + \ell(a) =
\ell(x) + \ell(a),$$ we deduce that $\ell(a) = 0$. Therefore $a =
e$ and we are done.

In the general case, let $x$ be a palindrome. There is
a central element $\Delta^{N} \in A_{\Gamma}^{+}$ such that
$\Delta^{2N}x \in A_{\Gamma}^{+}$. Since $\Delta^{2N}x$ is
still a palindrome, the previous argument applies. There is
$y \in A_{\Gamma}$ and $I \subseteq S$ such that
$\Delta^{2N}x = y \Delta_{I} \bar{y}$. Thus $x = \Delta^{-N}
y \Delta_{I} \bar{y} \Delta^{-N}$ is a desired
decomposition. \hfill $\blacksquare$\\

\subsubsection{Uniqueness} Let $x$ be a positive palindrome.
The set
$$ M(x) = \{ J \subseteq S \ | \ x = y \Delta_{J} \overline{y}\
{\rm{for}}\ {\rm{some}}\ y \in A_{\Gamma}^{+} \}$$ is
non-empty and finite. Since the ordering is total, there exists a
unique smallest element $\Delta_{I}$ such that $\Delta_{I} <
\Delta_{J}$ for all $J \in M(x)$, $J \not= I$. The subset
$$ N(x;I) = \{ y \in A_{\Gamma}^{+} \ | \ x = y \Delta_{I}
\overline{y} \}$$ is nonempty.  Since $<$ restricted to
$A_{\Gamma}^{+}$ is a well-ordering, $N(x;I)$ contains a unique
smallest element with respect to $<$.

Consider now a general palindrome $x \in A_{\Gamma}$. Then there
exists a central element $\Delta^{N} \in A_{\Gamma}^{+}$ such that
$\Delta^{2N} x \in A_{\Gamma}^{+}$. Clearly,
$\overline{\Delta^{2N} x} = \overline{x} \overline{\Delta^{2N}} =
x \Delta^{2N} = \Delta^{2N} x$. Applying the previous argument, we
obtain a canonical decomposition $\Delta^{2N} x = y \Delta_{I}
\overline{y}$ with $(\Delta_{I}, y)$ minimal among all other such
decompositions. Therefore, \begin{equation} x = \Delta^{-N} y
\Delta_{I} \bar{y} \Delta^{-N}. \label{eq:decomposition-g}
\end{equation}
Assume that there is another, distinct decomposition $x = z
\Delta_{J} \bar{z}$. Then we have $\Delta^{2N} x = \Delta^{N}z
\Delta_{J} \bar{z} \Delta^{N}  = y \Delta_{I} \bar{y}$. Since
$\Delta^{2N} x = y \Delta_{I} \bar{y}$ is the canonical
decomposition,  we have $(\Delta_{I}, y) \prec (\Delta_{J},
\Delta^{N}z)$. By left-invariance of $<$, this is equivalent to
$(\Delta_{I}, \Delta^{-N}y) \prec (\Delta_{J}, z)$. Hence, the
decomposition $(\ref{eq:decomposition-g})$ is unique. This is the
desired result.
\hfill $\blacksquare$\\

\begin{remark}
If $<$ does not restrict to a well-ordering on $A_{\Gamma}^{+}$,
we can still obtain uniqueness of the decomposition by requiring
the length (in Artin generators) of $\gamma$ to be minimal.
\end{remark}

\subsection{Image of the palindromization map}
\label{subsec:image}

Using Theorem \ref{th:dec-pal}, we must see that a palindrome $x =
\gamma \Delta_{I} \overline{\gamma}$ is pure if and only if
$\Delta_{I} = e$. The anti-automorphisms of $W_{\Gamma}$,
${\rm{rev}}:x \mapsto \overline{x}$ and $x \mapsto x^{-1}$
coincide on images of Artin generators in $W_{\Gamma}$. Hence they
coincide on $W_{\Gamma}$. Therefore, projecting $x$ to
$W_{\Gamma}$, we have $x = \gamma \Delta_{I} \gamma^{-1}$ (we
abusively keep the same notation for elements in $W_{\Gamma}$).
Now $x \in A_{\Gamma}$ is pure if and only if $\gamma \Delta_{I}
\gamma^{-1}$ is trivial in $W_{\Gamma}$. This occurs if and only
if $\Delta_{I}$ is trivial in $W_{\Gamma}$, hence trivial in
$A_{\Gamma}$ (by Tits' solution to the word problem).
Alternatively, one can verify inductively that for any $\varnothing \not= J
\subseteq S$, $\Delta_{J}$ is not pure. \hfill $\blacksquare$\\

\begin{remark}
The argument above yields a description of the images of
palindromes in $W_{\Gamma}$.
\end{remark}

\section{Applications} \label{sec:applications}

\begin{corollary}
Every element of order at most $2$ in $W_{\Gamma}$ lifts to a
palindrome in $A_{\Gamma}$.
\end{corollary}

\noindent{\bf{Proof}}. This is essentially a reformulation of the
previous observation (\S \ref{subsec:image}) coupled with the fact
that any element of order $2$ is the image of a conjugacy class of
$\Delta_{I}$ for
some subset $I$. \hfill $\blacksquare$\\

The following consequence of Theorem \ref{th:dec-pal} yields
restrictions on the possible fundamental elements occurring in the
canonical decomposition of a palindrome.

\begin{proposition} \label{prop:singleton}
Assume that $A_{\Gamma}$ is equipped with a left-invariant
ordering extending the subword order of $A_{\Gamma}^{+}$. Let $x$
be a palindrome in $A_{\Gamma}^{+}$
 and let $x = y \Delta_{I} {\overline{y}}$ be
its canonical decomposition. Then $I = \{s_1 \}  \cup \{ s_2 \}
\cup \ldots \cup \{ s_r \} \subseteq S$ and for all $1 \leq i < j
\leq r$, there is no edge in $\Gamma$ between $s_{i}$ and $s_{j}$.
In particular, $\Delta_{I} = \prod_{j} s_{j}$.
\end{proposition}

\noindent{\bf{Proof}}. Suppose that $I$ contains two non-commuting
Artin generators $s$ and $s'$. Denote by $m_{s,s'}$ the label of
the edge between $s$ and $s'$ in the Coxeter diagram. Then
$\Delta_{s, s'} = w_{m_{s,s'}}(s,s')$ divides $\Delta_{I}$. It
follows that $\Delta_{I} = s a s$ for some $s \in S$ and $a \in
A_{\Gamma}^{+}$. Since $\overline{\Delta_{I}} = \Delta_{I}$, we
have $\overline{a} = a$. Applying Theorem \ref{th:dec-pal} to $a$,
we obtain $a = b \Delta_{J} \bar{b}$ for some $b \in A_{\Gamma}$
and $J$ strictly contained in $I$. Hence $x = a' \Delta_{J}
\bar{a'}$ with $a' = s b \in A_{\Gamma}^{+}$ and $\Delta_{J}$
divides $\Delta_{I}$. Since the left-ordering of $A_{\Gamma}$
extends the subword order of $A_{\Gamma}^{+}$, we deduce that
$\Delta_{J} < \Delta_{I}$. This contradicts the minimality of the
canonical decomposition for $x$. \hfill $\blacksquare$

\begin{corollary} \label{cor:length}
Let $x = \gamma \Delta_{I} \overline{\gamma}$ be the canonical
decomposition of $x \in A_{\Gamma}$. Then
 $|I|$ is bounded by the number of the maximal subset of $S$ of commuting
 positive Artin generators.
\end{corollary}

\begin{example}
If $\Gamma = {\bf{A}}_{n}$, ${\bf{B}}_{n}$ or ${\bf{D}}_{n}$, then
$|I| \leq [\frac{n+1}{2}]$.
\end{example}

Recall that $\tau(x) = \Delta^{-1} x \Delta$, where $\Delta$ is
the fundamental element of $A_{\Gamma}$.


\begin{corollary} \label{cor:dec-special}
Let $x \in A_{\Gamma}$ be a palindrome. There is a decomposition
$x = y \Delta_{I} \bar{y}$ for some $I \subseteq S$, $y \in
A_{\Gamma}$ such that $\tau(\Delta_{I}) = \Delta_{I}$.
\end{corollary}

\noindent{\bf{Proof}}. It follows from \cite[\S 7]{BS}
that if $\tau$ is non trivial then each edge of the Coxeter
diagram is labelled by an odd integer.
Start with the canonical decomposition $x =
y \Delta_{I} \bar{y}$. By Proposition \ref{prop:singleton},
$\Delta_{I} = \prod_j s_j$ where $\{ s_j \}_j$ is a family of
commuting positive Artin generators. By definition of $\Delta$ as
a left and right least common multiple of $S$, the map $\tau: S
\to S, \ s \mapsto \tau(s)$ is a permutation of $S$ of order at
most $2$. Furthermore, since $\tau$ is a homomorphism, $\tau$ must
preserve the Coxeter diagram. Declare a subset $J \subseteq S$
admissible if there exists $z \in A_{\Gamma}$ such that $x = z
\Delta_{J} \bar{z}$. Let $J$ be an admissible subset such that
there is $s \in J$ commuting with all elements in $J$.
 Consider the
following operations on $J$.
\begin{enumerate}
\item[(A)] Adding a positive
Artin generator $s'$ to $J$ such that $s'$ commutes with
all elements of $J - \{ s\}$ and $s$ and $s'$ are joined by a
single edge in the Coxeter diagram.
\item[(B)] Replacing $s \in J$ by
another positive Artin generator $s' \in S$ such that $s'$
commutes with all elements of $J - \{ s \}$ and $s$ and $s'$ are
joined by a single edge in the Coxeter diagram.
\end{enumerate}
We claim that these operations do not affect admissibility. Denote
by $m_{s,s'} \in \Z$ the odd label of the edge between $s$ and $s'$.
Set $k = \frac{m_{s,s'}-1}{2}$.
We have
$$\begin{array}{rcl}
\Delta_{I \cup \{ s' \}} = \Delta_{(I - \{ s \}) \cup \{ s, s'\}}
= \Delta_{I - \{ s \}} \Delta_{s,s'} & = &
\Delta_{I - \{ s \}} w_{m_{s,s'}}(s,s') \\
& = & \Delta_{I - \{ s \}} w_{k}(s,s') s' \overline{w_{k}(s,s')} \\
& = & w_{k}(s,s') \Delta_{I - \{ s \}} s' \overline{w_{k}(s,s')} \\
& = & w_{k}(s,s') \Delta_{(I - \{ s \}) \cup s'}  \overline{w_{k}(s,s')}.
\end{array}$$
(We used the fact that $m_{s,s'}$ is odd in the fourth equality.) Thus
the operation (A) does not affect admissibility. The verification is
similar for operation (B). There is a sequence $(J_m)_{0 \leq m \leq k}$
of subsets of $S$ such that $J_0 = J$, each subset $J_{k+1}$ is obtained
from $J_k$ by means of one operation (A) or (B) and the final subset
$I=J_m$ satisfies $\tau(\Delta_{I}) = \Delta_{I}$. \hfill $\blacksquare$\\

\begin{theorem} \label{th:rev-tau-inv}
Let $x \in A_{\Gamma}$ such that ${\rm{rev}}(x) = x$ and $\tau(x)
= x$. Then there exists $I \subseteq S$, $y \in A_{\Gamma}$ such
that
\begin{equation}
 x = y \Delta_{I} \bar{y}, \ \ \tau(y) = y, \ \
\tau(\Delta_{I}) = \Delta_{I}. \label{eq:dec-tau-rev}
\end{equation}
\end{theorem}

\noindent{\bf{Proof}}. Suppose first that $x \in A_{\Gamma}^{+}$.
Lemma \ref{lem:prefix} yields $x = \Delta_{S(x)} y$ for some $y
\in A_{\Gamma}^{+}$. If $y = e$, we are done. Otherwise, $F(y)$ is
not empty and there is $s \in F(y)$. Since $\bar{x} = x$, $F(y)
\subseteq S(x)$. Hence $s \in S(x)$. Since $\tau(x) = x$,
$\tau(S(x)) = S(\tau(x)) = S(x)$. Hence $\tau(\Delta_{S(x)}) =
\Delta_{\tau(S(x))} = \Delta_{S(x)}$. It follows that $\tau(y) =
y$. We deduce that $s, \tau(s) \in F(y)$. Apply Lemma
\ref{lem:prefix} to $\bar{y}$: we have $x = \Delta_{s, \tau(s)} a
\Delta_{s, \tau(s)}$ for some $a \in A_{\Gamma}$. Since
$\Delta_{s, \tau(s)}$ is both rev- and $\tau$-invariant, we have
$\tau(a) = a = \bar{a}$. Furthermore, $\ell(a) < \ell(x)$. So we
can apply the argument again to $a$. This defines a recursive
procedure that stops if and only if the middle element $a$ either
trivial or is $\Delta_{I}$ for some subset $I \subseteq S$. This
is the desired result. For the general case, there is a central
element $\Delta^{N} \in A_{\Gamma}^{+}$ such that $z=\Delta^{2N} x
\in A_{\Gamma}^{+}$. Clearly $z$ is rev- and $\tau$-invariant. The
previous argument applies: there is $y \in A_{\Gamma}$ such that
$z = y \Delta_{I} \bar{y}$ for some $I \subseteq S$ and $\tau(y) =
y$ and $\tau(\Delta_{I}) = \Delta_{I}$. Hence $x = \Delta^{-N} y
\Delta_{I} \bar{y} \Delta^{-N} = y' \Delta_{I} \bar{y'}$ with $y'
= \Delta^{-N} y$ is a decomposition
with the required properties. \hfill $\blacksquare$\\

\begin{remark}
The two cases when the decomposition (\ref{eq:dec-tau-rev}) is
unique are $I = \varnothing$ and $I = S$, corresponding to Corollary
\ref{cor:pure-rev-tau}.
\end{remark}

\begin{remark} The subgroup $A_{\Gamma}^{\tau} \subset A_{\Gamma}$ of
$\tau$-invariant palindromes is again an Artin group of finite
type by a result due to J. Michel and P. Dehornoy -- L. Paris
\cite{DP}.
\end{remark}

\section{Further remarks} \label{sec:remarks}

\subsection{Left-orderings and palindromization}
Let $A_{\Gamma}$ be an Artin system equipped with a left-ordering
$<$ that has the property (SPPC). It is tempting to ask whether
the one-to-one palindromization map ${\mathfrak{pal}}:A_{\Gamma}
\to A_{\Gamma}, \ x \mapsto x \bar{x}$ is monotonic. It cannot be
decreasing since ${\mathfrak{pal}}(e) = e < x \bar{x} =
{\mathfrak{pal}}(x)$ for any $x > e$.

Below we show that for the classical braid group equipped
with Dehornoy ordering, the map ${\mathfrak{pal}}$ is {\emph{not}}
increasing. Set $x = s_{1} s_{2}$ and $y =
s_{1}^{2}$. It follows from the definition that $x < y$.
However, ${\mathfrak{pal}}(x) > {\mathfrak{pal}(y)}$. Indeed, we
have ${\mathfrak{pal}}(x) = s_1 s_{2}^{2} s_{1}$
and ${\mathfrak{pal}}(y) = s_{1}^{4}$. We rewrite $
{\mathfrak{pal}}(y)^{-1} {\mathfrak{pal}}(x)$ so as to find a
representative which is $s$-positive:
$$ \begin{array}{rcl}
{\mathfrak{pal}}(y)^{-1} {\mathfrak{pal}}(x) = s_{1}^{-4}
s_1 s_{2}^{2} s_{1} = s_{1}^{-3}
s_{2}^{2} s_{1} & = &
s_{1}^{-3} s_{2}^{2} s_{1} s_{2} s_{2}^{-1} \\
& = & s_{1}^{-3} s_2 s_1 s_2 s_1 s_{2}^{-1} \\
& = & s_{1}^{-3} s_1 s_2 s_1 s_1 s_{2}^{-1}\\
& = & s_{1}^{-2} s_{2} s_{1}^{2} s_{2}^{-1} \\
& = & s_{1}^{-2} s_{2} s_{1} s_{2} s_{2}^{-1}
s_{1} s_{2}^{-1} \\
& = & s_{1}^{-2} s_{1} s_{2} s_{1} s_{2}^{-1}
s_{1} s_{2}^{-1} \\
& = & s_{1}^{-1} s_{2} s_{1} s_{2}^{-1} s_{1}
s_{2}^{-1} \\
& = & s_{1}^{-1} s_{2} s_{1} s_{2} s_{2}^{-2}
s_{1} s_{2}^{-1} \\
& = & s_{1}^{-1} s_{1} s_{2} s_{1} s_{2}^{-2}
s_{1} s_{2}^{-1} \\
& = & s_{2} s_{1} s_{2}^{-2} s_{1} s_{2}^{-1}
\end{array}$$
is $s_{1}$-positive. Thus ${\mathfrak{pal}}(x) > {\mathfrak{pal}}(y)$.

\subsection{Artin groups of infinite type}

How much from the previous results remain true for Artin groups of
infinite type ? It is an open problem to determine which ones are
left-orderable (a fortiori to determine whether there is a
left-ordering that has (SPPC)).  On the other hand, Jacquemard's
algorithm is valid for all Artin groups. The key lemma of
\cite{DGKT} can be extended to all Artin monoids. Since all Artin
monoids embed naturally into their groups \cite{Paris}, the
palindromization map is one-to-one for all Artin groups.

\subsection{Garside groups}

Garside groups are a generalization of Artin groups of finite type
\cite{DP}. By means of the techniques used in this paper
(elementary divisibility theory and rewriting), the decomposition
for palindromes admits a generalization to Garside groups (no
uniqueness in general). However, the palindromization map is not
one-to-one in general. The simplest example is provided by the
Garside group
$$ G = \langle x, y\ | \ x^{2} = y^{2} \rangle. $$

\vskip 0.3cm

\noindent {\bf{Acknowledgements}}. We are indebted to Patrick
Dehornoy for his insight and for sharing his ideas about braid
ordering. We thank J. Birman, T. Fiedler and S. Orevkov for
valuable discussions.

\bibliographystyle{amsalpha}

%
%
%
%
%

\vskip 1cm

\end{document}